\documentclass[a4paper,12pt]{article}
\usepackage{amsmath}
\usepackage{amssymb}

\usepackage{natbib}
\bibpunct{[}{]}{;}{n}{,}{,}

\usepackage{amssymb, subfigure, epsfig, epsf, graphicx}

\usepackage{xspace}

\renewcommand{\leq}{\leqslant}

\newcommand{\ie}{i.e.\@\xspace}
\newcommand{\eg}{e.g.\@\xspace}

\begin{document}


\title{A quantitative investigation into the accumulation of rounding
      errors in the numerical solution of ODEs {\let\thefootnote\relax}
{\let\thefootnote\relax} } 

\date{\it University of Cambridge}

\author{Sebastian Mosbach, Amanda G. Turner}

\maketitle


\begin{abstract}
\noindent 
We examine numerical rounding errors of some deterministic solvers for systems of ordinary differential equations
(ODEs). We show that the accumulation of rounding errors results in a
solution that is inherently random and we obtain the theoretical
distribution of the trajectory as a function of time, the step
size and the numerical precision of the computer. We consider, in
particular, systems which amplify the effect of the rounding errors
so that over long time periods the solutions exhibit divergent behaviour. By
performing multiple repetitions with different values of the time step size, we observe numerically the random
distributions predicted theoretically. We mainly focus on the explicit
Euler and RK4 methods but also briefly consider more
complex algorithms such as the implicit solvers VODE and RADAU5.
\end{abstract}

{\it Abbreviated Title: Rounding errors in the numerical solution of ODEs}


\section{Introduction}
Consider ordinary differential equations (ODEs) of the form
$$
\dot{x}_t = b(x_t).
$$
These can be solved numerically using iteration methods of the type
$$
x_{t+h} = x_t + \beta(h,x_t)
$$
where $\beta(h,x) / h \to b(x)$ as $h \to 0$.

The simplest example of this is the Euler method where $\beta(h,x)=hb(x)$. This method is generally not used in
practice as it is relatively inaccurate and unstable compared to other methods. However more practical methods, such as
the fourth order Runge-Kutta formula (RK4), fall into this scheme.

When solving an ordinary differential equation numerically, each time
an iteration is performed, an error~$\epsilon$ is incurred due to rounding \ie
\begin{equation}
\label{introerr}
X^h_{t+h} = X^h_t + \beta(h, X^h_t) + \epsilon.
\end{equation}

Rounding errors in numerical computations are an inevitable
consequence of finite precision arithmetic. The first work 
thoroughly analysing the effects of rounding errors on numerical
algorithms is the classical textbook~\cite{Wil63}. A 
recent comprehensive treatment of the behaviour of numerical
algorithms in finite precision including an extensive list 
of references can be found in~\cite{Hig96}. Although rounding errors
are not random in 
the sense that the exact error incurred in any given calculation is
fully determined (see \cite{Hig96} or~\cite{For59}), probabilistic models
have been shown to adequately describe their behaviour.
In fact, statistical analysis of rounding errors can be traced back to one of
the first works on rounding error analysis~\cite{GolvNe51}. 

Henrici \cite{Hen62,Hen63,Hen64} proposes a probabilistic model for
the individual rounding errors, whereby they are independent and
uniform, the exact distribution depending on the specific finite
precision arithmetic being used. Using the central limit theorem he
shows that the theoretical distribution of the error accumulated
after a fixed number of steps in the numerical solution of an ODE is
asymptotically normal with variance proportional to $h^{-1}$. By
varying the initial conditions he obtains numerical distributions for
the accumulated errors, with good agreement. Hull and Swenson
\cite{HulSwe66} test the validity of the above model by adding a
randomly generated error with the same distribution at each stage of
the calculation, and 
comparing the distribution of the accumulated errors with those
obtained purely by rounding. They observe that although rounding is
neither a random process nor are successive errors independent,
probabilistic models appear to provide a good description of what
actually happens. 

We shall concentrate on floating point arithmetic, as used by modern
computers. However, our methods can be used equally well for any
finite precision arithmetic. We use the model, discussed and tested by
the authors above, whereby under generic conditions, the errors in
\eqref{introerr} can be viewed as independent, zero mean, uniform random
variables $\epsilon_i \sim U[-|X^h_{t,i}| 2^{-p}, |X^h_{t,i}|
 2^{-p}]$, $p$ being a constant determined by the precision of the computer.

In the first half of the paper we analyse the cumulative effect of
these rounding errors as the step size $h\to 0$. Where previous
authors have considered the accumulated error at a particular point,
we derive a theoretical model for the entire trajectory. Cases in
$\mathbb{R}^2$ where the ordinary differential equation has a saddle
fixed point at the origin exhibit the 
most interesting behaviour, as the structure of the ODE system
amplifies the effect of the rounding errors and causes the numerical
solution to diverge from the actual solution. We show
that in this case the solution $X^h_t$ is inherently random and we obtain its theoretical distribution as an explicit
function of time, the step size and the precision of the computer. We shall see that as the step size $h\to 0$, the
numerical solution exhibits three types of behaviour, depending on the time. More precisely, there exists a
constant~$c$, determined by the ODE system, such that for times much smaller than $- c \log h$ the numerical solution
converges to the actual solution; for times close to $- c \log h$ the solution undergoes a transition, driven by a
Gaussian random variable whose distribution we shall obtain; for times
much larger than $- c \log h$ the numerical solution
diverges from the actual solution.

In the second half of the paper, we perform numerical simulations which illustrate this behaviour. By performing
multiple repetitions with different values of the time step size, we observe the random distributions predicted
theoretically. Where previous authors have obtained their numerical
distributions by varying the initial conditions, we do so by introducing small
variations in the step size $h$.  We show that during the transition
period described in the previous paragraph, the numerical solution
intersects straight 
lines through the origin and we compare the theoretical and numerical distributions for the points at which these
intersections occur. Both the mean and the standard deviation of these distributions are of the form $a h^{\gamma}$,
where $ \gamma \in (0, 1/2]$ is a constant determined by the ODE system, and~$a$ can be found explicitly in terms of
the precision of the computer, \ie the number of bits used internally by the computer to represent floating point
numbers. We mainly focus on the explicit Euler and RK4 methods, but
show that the same behaviour is also observable for more 
complex algorithms such as the implicit solvers VODE and RADAU5.

\section{Theoretical background}
In the paper \cite{Tur04}, limiting results are established for sequences of Markov processes that approximate solutions of
ordinary differential equations with saddle fixed points. We shall outline these results and then show that the
rounding errors accumulated when performing numerical schemes for solving ordinary differential equations can be viewed
as a special case of this. This enables us to quantify how the rounding errors combine, and show that the resulting
numerical solutions exhibit random behaviour, the exact distribution of which is obtained.

\subsection{Behaviour of stochastic jump processes}

We are interested in ordinary differential equations of the form
\begin{equation}
 \label{ode}
 \dot{x}_t = b(x_t).
\end{equation}
We focus on $\mathbb{R}^2$ in the case where the origin is a saddle fixed point of the system \ie $b(x_t) = B
x_t + \tau(x_t)$, where $B$ is a matrix with eigenvalues $\lambda, -\mu$, with $\lambda, \mu > 0$ and $\tau(x) =
O(|x|^2)$ is twice continuously differentiable. This case is of particular interest as the structure of the system
amplifies the effect of the rounding errors and causes the numerical
solution to diverge from the actual solution over large times. Similar
behaviour can be observed in higher dimensions where the matrix $B$
has at least one positive and one negative eigenvalue, although the
corresponding quantitative analysis is much harder and we do not go into it here.

\begin{figure}[h]
 \centering
     \epsfig{file=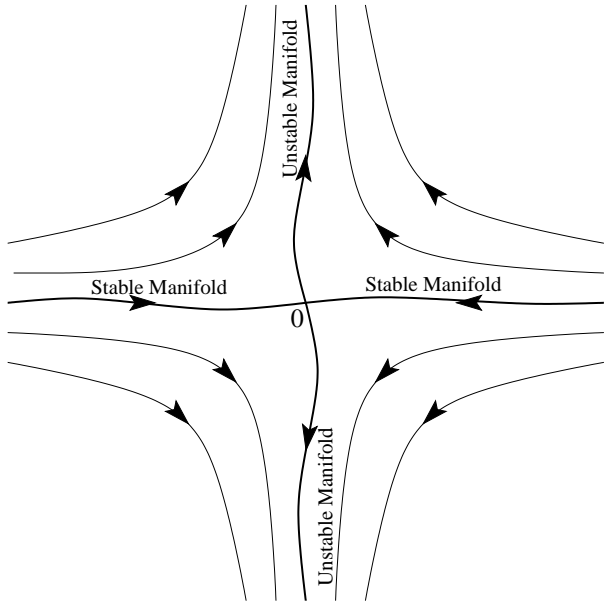, width=8cm}
  \caption{\textsl{The phase portrait of an ordinary differential equation having a saddle fixed point at the origin (taken from~\cite{Tur04}).}}
  \label{phaseok}
\end{figure}

The phase portrait of \eqref{ode} in the neighbourhood of the origin
is shown in Figure \ref{phaseok}.
In particular, there exists some $x_0 \neq 0$ such that $\phi_t(x_0) \rightarrow 0$ as $t \rightarrow \infty$, where
$\phi$ is the flow associated with the ordinary differential equation \eqref{ode}. The set of such $x_0$ is the stable
manifold.  There also exists some $x_{\infty}$ such that $\phi^{-1}_t(x_{\infty}) \rightarrow 0$ as $t \rightarrow
\infty$. The set of such $x_{\infty}$ is the unstable manifold.

\begin{figure}[h]
 \centering
     \epsfig{file=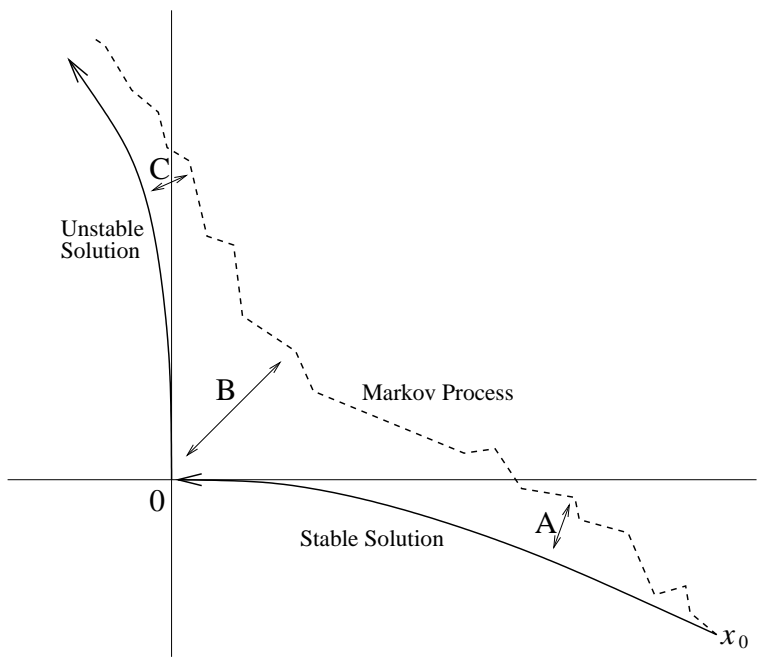, width=10cm}
  \caption{\textsl{Diagram showing how the Markov process $X_t^N$ deviates from the stable solution $\phi_t(x_0)$ for
    large values of $t$ (taken from~\cite{Tur04}).}}
 \label{over}
\end{figure}

Fix an $x_0$ in the stable manifold and consider sequences $X_t^N$ of Markov processes starting from $x_0$, which
converge to the solution of \eqref{ode} over compact time intervals. The processes are indexed so that the variance of
the fluctuations of $X_t^N$ is inversely proportional to $N$. If we allow the value of $t$ to grow with $N$ as a
constant times $\log N$, $X_t^N$ deviates from the stable solution to a limit which is inherently random, before
converging to an unstable solution (see Figure \ref{over}).
More precisely, we observe three different types of behaviour depending on the time scale:
\begin{enumerate}
 \item[A.]
  On compact time intervals, $X_t^N$ converges to the stable solution
  of \eqref{ode}, the fluctuations around this limit being of order
  $N^{-\frac{1}{2}}$. The exact distribution of the fluctuations is
  asymptotically $N^{-\frac{1}{2}} \gamma_t$ where $\gamma_t$ is the
  solution to a linear stochastic differential equation, described in
  \cite{Tur04}. 
 \item[B.]
  Let $v_1$ and $v_2$ be the unit eigenvectors of $B$ corresponding to
  $- \mu$ and $\lambda$ respectively.
  There exists some $\overline{x}_0 \neq 0$, depending only on $x_0$,
  and a Gaussian random variable $Z_{\infty}$, such that if $t$ lies
  in the interval $[R, \frac{1}{2 \lambda} \log N - R]$, then
  $$
  X_t^N = \overline{x}_0 e^{- \mu t} (v_1 + \epsilon_1) +
  N^{-\frac{1}{2}}Z_{\infty} e^{\lambda t}(v_2 + \epsilon_2)
  $$
  where $\epsilon_i(t,N) \rightarrow
  0$ uniformly in $t$ in probability as $R, N \rightarrow \infty$. In
  other words, $X_t^N$ can be approximated by the solution to the
  linear ordinary differential equation
  \begin{equation}
   \label{linode}
   \dot{y_t} = B y_t
  \end{equation}
  starting from the random point
  $\overline{x}_0 v_1 + N^{-\frac{1}{2}}Z_{\infty} v_2$.
 \item[C.]
  Provided $Z_{\infty} \neq 0$, on time intervals of a fixed length
  around $\frac{1}{2 \lambda} \log N$,
  $X_t^N$ converges to one of the two unstable solutions of
  \eqref{ode}, each with probability~$1/2$, depending on the sign of
  $Z_{\infty}$. 
\end{enumerate}

\subsection{Accumulation of rounding errors}
\label{accrou}
We can apply the above results to describe quantitatively how rounding errors accumulate when solving ordinary
differential equations of the form \eqref{ode} numerically. In particular we consider using iteration methods of the
type
$$
x_{t+h} = x_t + \beta(h,x_t)
$$
where $\beta(h,x) / h \rightarrow  b(x)$ as $h \rightarrow 0$ \eg the Euler method where $\beta(h,x)=hb(x)$.

Each time an iteration is performed, an error $\epsilon = \epsilon(h,t)$ is incurred due to rounding, so we obtain a
process $(X_t^h)_{t \in h\mathbb{N}}$ iteratively by
\begin{equation}
\label{erroracc}
 X^h_{t+h} = X^h_t + \beta(h, X^h_t) + \epsilon.
\end{equation}
Modern computers store real numbers by expressing them in binary as $x = m 2^n$ for some $1 \leq |m| < 2$ and $n \in
\mathbb{Z}$, and allocate a fixed number of bits to store the
mantissa $m$ and a (different) fixed number of bits to store the
exponent $n$ \cite{IEEE754}. When adding a smaller number to $x$, the size of the rounding error incurred is between 0
and $2^n.2^{-p} =2^{\lfloor \log_2 |x| \rfloor - p}$, where $p$ is the number of bits allocated to store the mantissa.
Although it is possible to carry out the calculations below using the exact value of $2^{\lfloor \log_2 |x| \rfloor -
p}$, the calculations are greatly simplified by approximating it by $|x| 2^{-p}$. This results in the `effective' value
of $p$ differing from the actual value of $p$ by some number between 0 and 1. Under generic conditions, the errors
$\epsilon$ can therefore be viewed as independent, mean zero, uniform random variables with approximate distribution
$\epsilon_i \sim U[- | X^h_{t,i}| 2^{-p}, |X^h_{t,i}| 2^{-p}]$ (see \cite{Hen62,Hen63,Hen64}). The assumption that the $\epsilon_i$ are independent is
violated in certain pathological cases, for example where there is a lot of symmetry in the components. However, in
general it is a reasonable assumption.

Although the above iterations are carried out at discrete time intervals, it is convenient to embed the processes in
continuous time by performing the iterations at times of a Poisson process with rate $h^{-1}$. As $\beta(h,x)$ does not
depend on $t$, this does not affect the shape of the resulting trajectories. In this way we obtain Markov processes
$X_t^h$ that approximate the stable solution of \eqref{ode} for small values of $h$. If, in addition, we assume that
$$
h^{-\frac{1}{2}}\left ( \frac{\beta(h,x)}{h} - b(x) \right ) \rightarrow 0
$$
as $h \rightarrow 0$ (note that both Euler and Runge-Kutta satisfy this condition), then under the correspondence $N
\sim h^{-1}$, we satisfy the conditions needed to apply the results in \cite{Tur04}. Our numerical solution therefore
exhibits the following random behaviour:
\begin{enumerate}
 \item[A.]
  For times of order much smaller than $- \log h$, $X_t^h$
  approximates the stable solution of \eqref{ode}, the fluctuations
  around this limit being of order $h^{\frac{1}{2}}$.
 \item[B.]
  There exists some $\overline{x}_0 \neq 0$, depending only on $x_0$,
  and a Gaussian random variable $Z_{\infty}$, such that if $t$ lies
  in the interval $[-c \log h, -\frac{1}{2 \lambda} \log h + c \log
  h]$ for some $c > 0$, then $X_t^h$ is asymptotic to
  \begin{equation}
  \label{randlin}
   \overline{x}_0 e^{- \mu t}v_1 + h^{\frac{1}{2}}Z_{\infty} e^{\lambda t}v_2,
  \end{equation}
  the solution to the linear ordinary
  differential equation \eqref{linode} starting from the random point
  $\overline{x}_0 v_1 + h^{\frac{1}{2}}Z_{\infty} v_2$.
 \item[C.]
  Provided $Z_{\infty} \neq 0$, on time intervals around $- \frac{1}{2
  \lambda} \log h$ whose length is of much smaller order than $- \log
  h$, $X_t^h$ approximates one of the two unstable solutions of
  \eqref{ode}, each with probability $\frac{1}{2}$, depending on the
  sign of $Z_{\infty}$.
\end{enumerate}

The random behaviour resulting from the accumulation of rounding errors is most noticeable on time intervals of fixed
lengths around $- \frac{1}{2(\lambda + \mu)} \log h$, as for these values of $t$ the two terms $\overline{x}_0 e^{- \mu
t}$ and $h^{\frac{1}{2}}Z_{\infty} e^{\lambda t}$ in \eqref{randlin} are of the same order. Over this time interval,
the numerical solution undergoes a transition from converging to the
actual solution to diverging from it. During this transition, for each
value of $\theta \in (0, \pi/2)$, $X_t^h$ crosses one of the straight
  lines passing through 0 in the direction $\cos \theta v_1 \pm \sin \theta
v_2$. These intersections are important as 
they indicate the onset of divergent behaviour. The distribution of
the point at which $X_t^h$ intersects one of the lines in the
direction $\cos \theta v_1 \pm \sin \theta v_2$ is asymptotic to 
\begin{equation}
\label{intersect}
 h^{\frac{\mu}{2(\lambda + \mu)}}|Z_{\infty}|^{\frac{\mu}{\lambda +
 \mu}}|\overline{x}_0|^{\frac{\lambda}{\lambda + \mu}} |\tan
 \theta|^{\frac{\mu}{\lambda + \mu}} (\cos \theta v_1 \pm \sin \theta v_2).
\end{equation}

In Section \ref{zinfcalc} we show how to evaluate the variance of $Z_{\infty}$, doing so explicitly in the linear case
and obtaining bounds in the non-linear case. In Section \ref{NumExp} we verify these results by numerically obtaining
the predicted distribution for hitting a line through the origin.

\subsection{Explicit calculation of the variance}
\label{zinfcalc}

Suppose that we are using a numerical scheme that satisfies the above conditions to obtain a solution to the ordinary
differential equation \eqref{ode} starting from $x_0$ for some $x_0$ in the stable manifold. In the non-linear case we
require that $x_0$ is sufficiently close to the origin such that $\tau(x_0)$ is small. In general, for simplicity, we
shall assume that $|x_0| \leq 1$.

We define the flow $\phi$ associated with this system by
$$
\dot{\phi}_t(x) = b(\phi_t(x)), \quad \phi_0(x) = 0
$$
and let $x_t = \phi_t(x_0)$.

Suppose $v_1$, $v_2 \in \mathbb{R}^2$ are the unit right-eigenvectors of $B$ corresponding to $- \mu$, $\lambda$
respectively, and that $v'_1, v'_2 \in (\mathbb{R}^2)^{\ast}$ are the corresponding left-eigenvectors (i.e. $v_i' v_j =
\delta_{ij}$).

Define
$$
\overline{x}_0 = \lim_{t \rightarrow \infty} e^{\mu t} v'_1 \phi_t(x_0)
$$
and
$$
D_s = \lim_{t \rightarrow \infty} e^{- \lambda t} v'_2 \nabla \phi_t (x_s).
$$
It is shown in \cite{Tur04} that these limits exist and that $|\overline{x}_0| \leq 2|x_0| \leq 2$, and $|D_s| \leq 2$.

Finally, let
$$
a(x) = \frac{1}{3} 2^{-2p} \left ( \begin{matrix} x_1^2 & 0 \\ 0 &
       x_2^2 \end{matrix} \right )
$$
be the covariance matrix of the multivariate uniform random variable $\epsilon$, defined in equation \eqref{erroracc},
when $X_t^h=x$. Then $Z_{\infty} \sim N (0, \sigma_{\infty}^2)$, where
$$
\sigma^2_{\infty} = \int_0^{\infty} e^{- 2 \lambda s} D_s a(x_s) D_s^{\ast} ds.
$$
Note that $\sigma^2_{\infty} \leq \frac{2}{3 \lambda}2^{-2p}$.

In the general non-linear case, evaluating $\sigma_{\infty}^2$ explicitly is not possible as it involves solving
\eqref{ode}. It is possible to obtain a better approximation than that above, although the important observation is
that it is proportional to $2^{-2p}$.

In the linear case, $\phi_t(x) = e^{B t} x$ and $x_0 = |x_0| v_1$. Hence $x_t = |x_0| e^{- \mu t} v_1$, $\overline{x}_0
= |x_0|$, and $D_s = v_2'$, and so
$$
 \sigma_{\infty}^2 = \frac{1}{3(\lambda + \mu)}2^{-2p} |x_0|^2
 (v_{1,1}v'_{2,1})^2.
$$
Note that the directions of $v_1$ and $v_2'$, relative to the standard
basis, are critical. For example, if either $v_1$ or
$v_2'$ is parallel to one of the standard basis vectors, then $\sigma_{\infty}^2=0$.

\section{Numerical experiments}
\label{NumExp}
In this section we solve ODEs numerically using deterministic solvers and observe the predicted random distributions
arising as a consequence of the accumulation of rounding errors. For simplicity, and in order to observe the desired
effects as clearly as possible, we mainly focus on the most elementary of all numerical ODE solution methods: the
standard explicit Euler algorithm with constant time step
size. However, we observe similar behaviour for RK4 and also briefly mention
results obtained with 
more complex solvers, such as VODE \cite{BroByrHin89}. For a recent overview of ODE solvers see \cite{Cas03}; for an
introductory text see \cite{HaiWan96}.
\subsection{The system}

\begin{figure}[h]
  \centering
  \epsfig{file=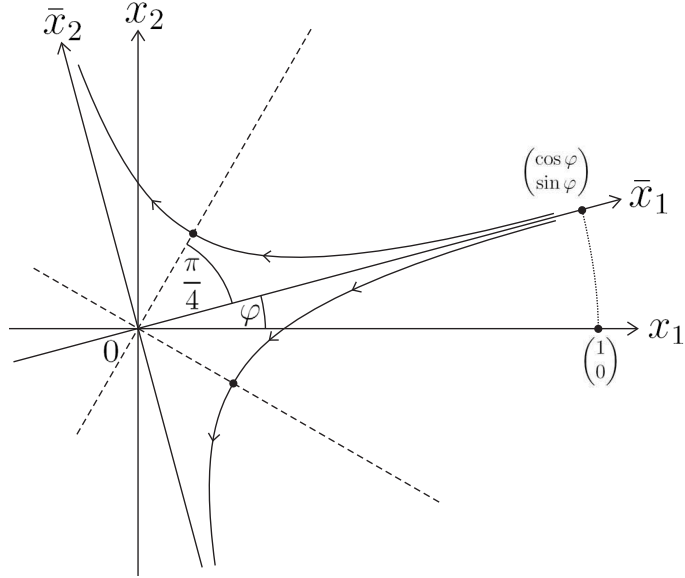,width=9.0cm}
  \caption{\textsl{Phase space for the saddlepoint ODE system
    \eqref{EqnODE} with sample trajectories and lines where hitting
    distributions are recorded (dashed lines).}}
  \label{FigSaddlepoint}
\end{figure}

For $x:[0,\infty) \rightarrow \mathbb{R}^2$, consider the linear ODE
\[\dot{x}(t)=Bx(t)\]
where \[B=\begin{pmatrix}
  -\mu & 0\\
  0 & \lambda
  \end{pmatrix}\] for fixed $\lambda,\mu>0$. We introduce new coordinates
\[\bar{x}(t)=R(\varphi)x(t)\]
by rotating about the origin by a fixed angle $\varphi\in[0,\pi/2)$, \ie
\[R(\varphi)=\begin{pmatrix}
  \cos\varphi & -\sin\varphi\\
  \sin\varphi & \cos\varphi
  \end{pmatrix}.\]
We arrive at the transformed system
\begin{equation}
  \label{EqnODE}
  \dot{\bar{x}}(t)=\bar{B}(\varphi)\bar{x}(t)
\end{equation}
with
\[\bar{B}(\varphi)=R(\varphi)BR(\varphi)^\top,\]
which will be the system under consideration in the following. Throughout, we use as initial value
\begin{equation}
  \label{EqnInitialValue}
  \bar{x}(0)=R(\varphi)\begin{pmatrix}1\\0\end{pmatrix}=\begin{pmatrix}\cos\varphi\\\sin\varphi\end{pmatrix}.
\end{equation}
The phase space evolution is sketched in Figure \ref{FigSaddlepoint}.

\subsection{Theoretical hitting distribution}
\label{SubSubSecThHitDistr}
As discussed in Section \ref{accrou}, the numerical solution to the
above ODE system undergoes a transition from converging to the actual
solution to diverging from it. During this transition, the numerical
trajectory crosses one of the straight lines passing
through 0 at an angle $\phi \pm \theta$ for each value of $\theta \in
(0, \pi/2)$. These 
intersections are important as they indicate the onset of
divergent behaviour. The hitting distributions also provide a means of 
measuring the random variable $Z_{\infty}$, which drives the random
variations in our solutions, and hence of verifying the theoretical
results. 

Equation \eqref{intersect} gives the asymptotic distribution of the
magnitude of the point at which the numerical solution
hits the line through the origin at angle $\varphi \pm \frac{\pi}{4}$ as
$|Z|^{\frac{\mu}{\lambda + \mu}}$ where $Z$ is a Gaussian random variable with mean 0 and variance
\begin{equation}
\label{sigmathe}
 \sigma^2 = h \sigma_{\infty}^2 = \frac{1}{3(\lambda + \mu)} h 2^{-2p}
 (\cos \varphi \sin \varphi)^2  
\end{equation}
\ie $Z \sim\mathcal{N}(0,\sigma^2)$. We obtain an explicit formula for the asymptotic distribution by starting from the
$\mathcal{N}(0,\sigma^2)$ distribution
\[p(x)\text{d}x=\frac{1}{\sqrt{2\pi} \sigma}\exp\Big(-\frac{1}{2
  \sigma^2}x^2\Big)\text{d}x\]
and performing a transformation of the variable given by $y=|x|^{\frac{\mu}{\lambda + \mu}}$. The result is:
$$
  p(y)\text{d}y=\frac{2(\lambda + \mu)}{\sqrt{2 \pi} \sigma
  \mu}y^{\frac{\lambda}{\mu}}\exp\Big(-\frac{1}{2
  \sigma^2}y^{\frac{2 (\lambda + \mu)}{\mu}}\Big)\text{d}y.
$$
In the case $\lambda = \mu = 1$ which we shall consider below, setting $a = \frac{4}{\sqrt{2 \pi} \sigma}$, we obtain
the family of distributions
\begin{equation}
  \label{EqnFitDistr}
  f(x) \text{d}x = a x \exp\Big(-\frac{\pi}{16} a^2 x^{4}\Big)\text{d}x
\end{equation}
which we shall fit to our numerical data to confirm our theoretical value of $a$.

\subsubsection{Choice of parameters}
Since there is no random element contained in a deterministic solver, for each repetition at least one parameter has to
vary. At first, we tried varying the initial value in the direction of the eigenvectors, but this did not yield any
interesting results. The chosen distribution of initial values was reproduced exactly in the hitting distribution and
no randomness could be observed. Also, from an aesthetic point of view, it is preferable to vary only `internal', \ie
numerical parameters of an algorithm, such as the time step size or error tolerances, instead of varying `physical'
parameters of the system like the initial value.

The only internal parameter of Euler's algorithm is the time step size~$h$, which we varied as follows. Given a
user-supplied value of~$h$, the step size~$h_i$ for the~$i^\text{th}$ repetition ($i\in\{1,\dots,L\}$) is defined by
\[h_i=h+\Delta h(i-1-k),\]
where the number of repetitions $L=2k+1$ and $0 < \Delta h \ll h$ are also user-supplied. For all simulations, we set
$k=10^4$.

In each run, the trajectory at some point intersects the lines $\bar{x}_1=\pm\bar{x}_2$ (the dashed lines in Figure
\ref{FigSaddlepoint}). In order to produce histograms, we partition the interval $[0,1]$ into a given fixed number of
subintervals of equal length and count how many times~$y$ falls into each subinterval, where~$y$ denotes the distance
of the point of intersection from the origin.

Since the limit distribution is given by $|Z|^{\frac{\mu}{\lambda+\mu}}$, if the values of~$\lambda$ and~$\mu$ differed
significantly then the distribution would be hard to observe in a numerical experiment. Furthermore, if $\lambda\gg\mu$
then the trajectories are very quickly pushed away from the $\bar{x}_1$-axis so that fluctuations (\ie rounding errors)
become largely irrelevant. Conversely, if $\mu\gg\lambda$ then the amplification of deviations from the
$\bar{x}_1$-axis is too weak to be observable. These arguments suggest choosing $\lambda$ and~$\mu$ within the same
order of magnitude, and we therefore chose $\lambda=\mu=1$ for all simulations.

Another subtlety concerns the choice of the rotation angle~$\varphi$. For certain values, trivial trajectories or
symmetry effects can occur which conceal the desired accumulation of rounding errors. For instance, for~$\varphi=0$ the
second component~$\bar{x}_2$ of the solution is always zero, and therefore the trajectory stays on the
line~$\bar{x}_2=0$ (or equivalently~$x_2=0$) with no fluctuations. Note that this is in agreement with
$\sigma^2=0$ in equation \eqref{sigmathe}. For $\varphi=\pi/4$, any rounding error that appears in one
component also appears in the other one, which implies that, again, the trajectory always stays on the
line~$\bar{x}_2=0$ (or equivalently~$x_1=x_2$). This case is pathological as it consistently violates our assumption
that the rounding errors are independent for each component. For these reasons, we chose $\varphi=\pi/5$ throughout.

Reasonable choices of~$h$ and~$\Delta h$ are limited by several factors. If~$h$ is too large (in the considered case
$h>10^{-1}$ for both single and double precision) then the observed hitting distributions differ substantially from the
theoretical one because not enough rounding errors can accumulate and hence the deviations are not sufficiently random.
The onset of such effects can be seen for large values of~$h$ in Figure \ref{FigsAStepSize}. Lower bounds on~$h$ are
imposed on the one hand by computational cost and on the other hand by the numerical accuracy of the computer on which
the calculations are performed. In practice, however, computational expense becomes prohibitive for values of~$h$ much
larger than the smallest values permitted by numerical accuracy. Our particular choice of step size distribution
requires that~$k\Delta h$ should be (much) smaller than~$h$. The lower limit for~$\Delta h$ is determined solely by the
numerical precision, \ie $\Delta h/h$ must not be smaller than the numerical precision because otherwise there is no
variation at all and hence no randomness.

It is beyond the scope of this paper to investigate in detail the dependence of our observations on the distribution of
step sizes. However, preliminary experiments with varying~$\Delta h$ and even with non-uniform step size distributions
suggest that this dependence is very weak for a wide range of conditions. The mean value of the step size distribution,
whose effect on the shape of the hitting distribution (\ie the parameter~$a$ in equation \eqref{EqnFitDistr}) has been
demonstrated in Figure \ref{FigsAStepSize}, seems to be most significant. On the other hand, the variance seems
irrelevant apart from being large enough to produce randomness, and does not affect the parameter~$a$. This is also
supported by the fact that~$a$ is (asymptotically) independent of the
number of repetitions~$L$. 

\begin{figure}[h]
  \centering
  \epsfig{file=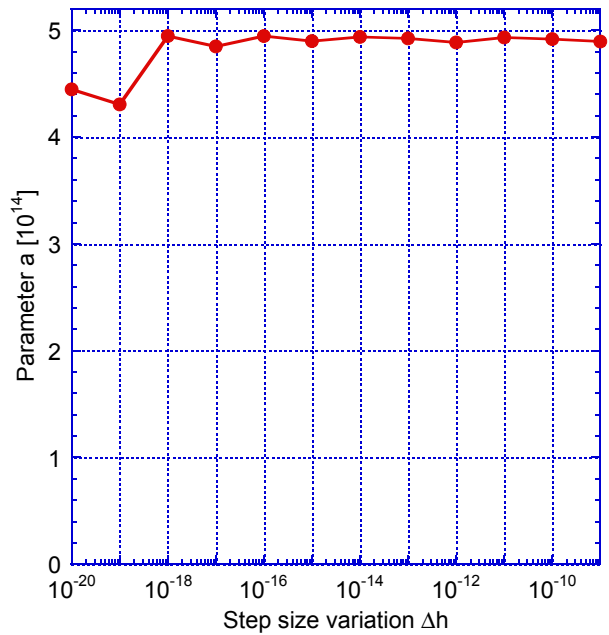,width=8.0cm}
  \caption{\textsl{Step size variation for Euler's algorithm (double
  precision, step size~$h=10^{-4}$, $L=20001$ repetitions each).}}
  \label{FigEulerDhVar}
\end{figure}

Figure \ref{FigEulerDhVar} shows that the shape of the distribution exhibits no discernible systematic dependence
on~$\Delta h$ over at least nine orders of magnitude. This reinforces that the parameter~$a$ is unaffected by the
variance of the step size distribution for this particular case. The deviations seen for values of~$\Delta h$ smaller
than about $10^{-19}$ are due to the fact that~$\Delta h/h$ approaches
the limits of numerical precision.

\subsubsection{Results and observations for explicit methods}

\begin{figure}[h]
  \centering
  \epsfig{file=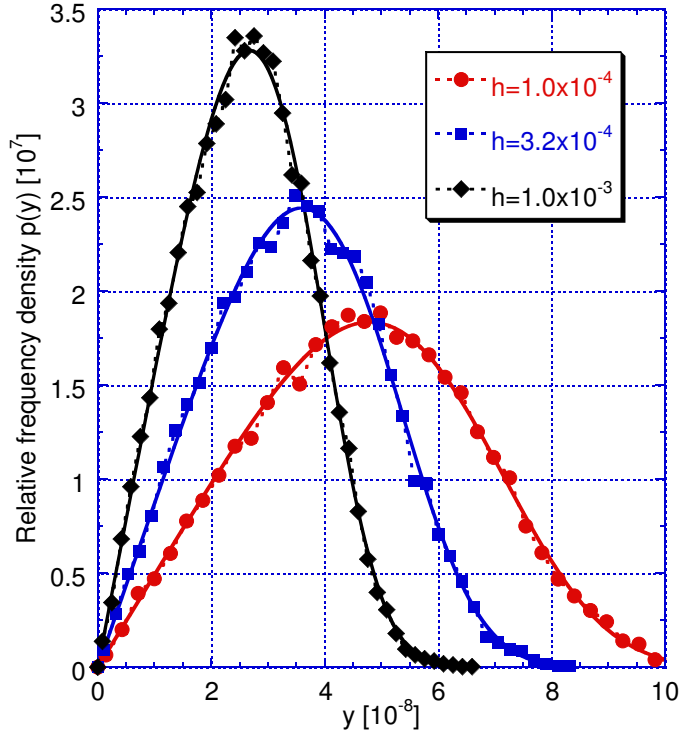,width=9.0cm}
  \caption{\textsl{Observed hitting distributions with theoretical fits for
Euler's algorithm ($\Delta h=10^{-10}$, $L=20001$ repetitions each).}}
  \label{FigEulerHitDistr}
\end{figure}

Using the given values, distributions as shown in Figure \ref{FigEulerHitDistr} are obtained. We fitted the theoretical
distribution \eqref{EqnFitDistr} to the ones produced numerically with very good agreement.
%
\begin{figure}[h]
  \centering
  \mbox{
    \subfigure[Single precision ($\Delta h=10^{-8}$).]
      {\label{FigAStepSizeSgl}
      \epsfig{file=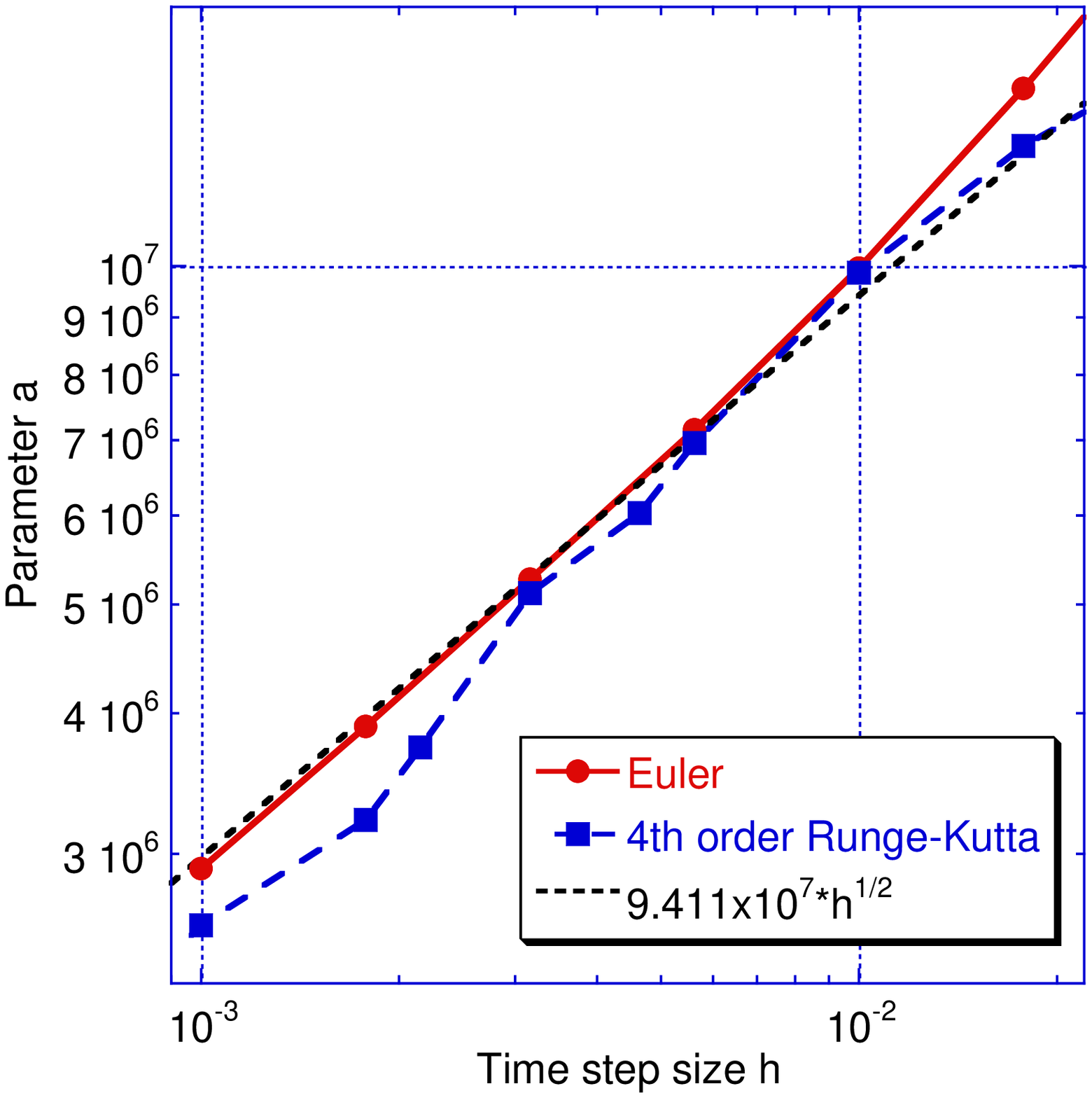,width=7.12cm}}
    \subfigure[Double precision ($\Delta h=10^{-10}$).]
      {\label{FigAStepSizeDbl}
      \epsfig{file=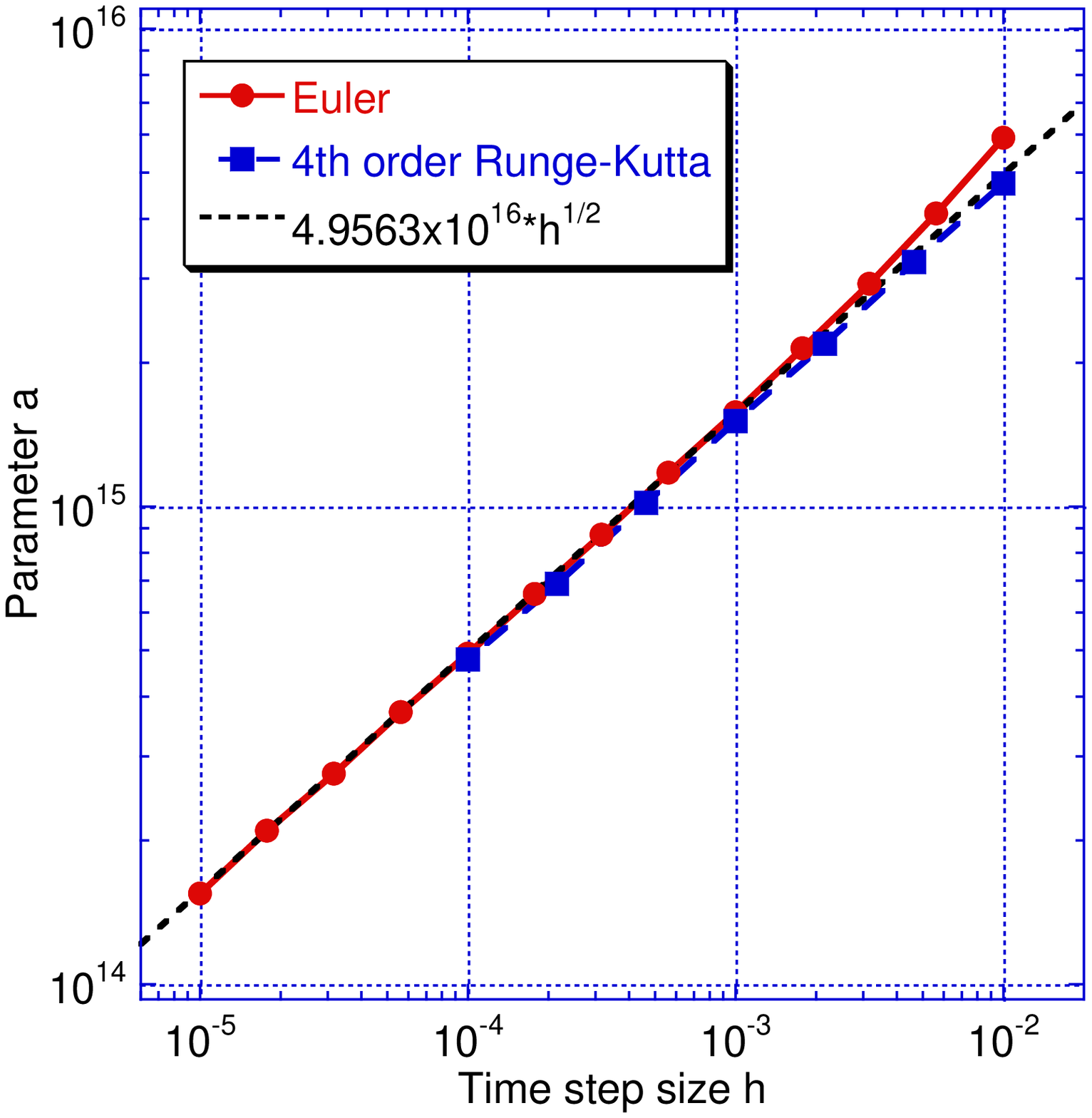,width=7.0cm}}}
  \caption{\textsl{Parameter~$a$ in equation \eqref{EqnFitDistr} as function of the time step size~$h$ for simple
    explicit methods (Euler and $4^\text{th}$ order Runge-Kutta).}}
  \label{FigsAStepSize}
\end{figure}

From the numerical experiments we obtain for each~$h$ a distribution
of the form \eqref{EqnFitDistr}, where the
parameter~$a$ is given as a result of the fitting procedure. In Figure \ref{FigsAStepSize} the parameter~$a$ is plotted
as a function of the time step size~$h$, both for single (Figure \ref{FigAStepSizeSgl}) and double (Figure
\ref{FigAStepSizeDbl}) precision ($4$ and $8$ bytes internal representation of floating point numbers respectively).
Error bars due to the fit are only about~$1\%$ and hence invisibly small. In both cases, the dependence between~$a$
and~$h$ seems to be well described by $a\propto\sqrt{h}$. Intuitively, one might explain the qualitative behaviour as
follows: The smaller the step size, the more numerical errors accumulate, and hence the broader the distribution. The intuition ``the smaller the step size, the higher the accuracy, and hence the narrower the distribution'',
although at first sight equally valid, is false. In the single precision case, because of the lower accuracy
compared to double precision, the distributions are much broader for a given step size.

Equation \eqref{sigmathe} predicts the value of $ah^{-\frac{1}{2}}$ to be
$$
ah^{-\frac{1}{2}} = \frac{4 \sqrt{3}}{\sqrt{\pi} \cos \frac{\pi}{5} \sin
  \frac{\pi}{5}}\times 2^p = 8.220 \times 2^p.
$$
For Euler's method, the above data give $ah^{-\frac{1}{2}} =9.411 \times 10^7$ for single precision and $ah^{-\frac{1}{2}} = 4.956 \times
10^{16}$ for double precision. For $4^\text{th}$ order Runge-Kutta,
the values are $ah^{-\frac{1}{2}} = 9.27 \times 10^7$ (with a
relatively large error of $\pm 0.12 \times 10^7$) for single precision
and $ah^{-\frac{1}{2}} = 4.746 \times 10^{16}$ for double
precision. Using the approximation discussed in Section \ref{accrou},
the actual value of $p$ is between 23 and 24, when working in single
precision, and between 52 and 53 when working in double precision, the 
particular value depending on the exact number being computed. Our
theoretical results therefore predict $a 
h^{-\frac{1}{2}}$ lies between $6.895 \times 10^7$ and $1.379 \times 10^8$ for single precision and between $3.702
\times 10^{16}$ and $7.404 \times 10^{16}$ for double precision.

There are three possible sources of error in our calculations. The first is the error in fitting the numerical data to
the theoretical model, the second is that our theoretical models are based on asymptotic results as $h \rightarrow 0$,
whereas we are applying them to values of $h$ which are necessarily larger than the precision of the computer. The
third source of error arises from the assumption that at each stage the rounding error can be viewed as an independent
uniform random variable, depending on a fixed value of $p$. The above results show the above errors are all small and
our theoretical model provides a very good fit.

\subsubsection{Implicit solvers}
%
\begin{figure}[h]
  \centering
  \epsfig{file=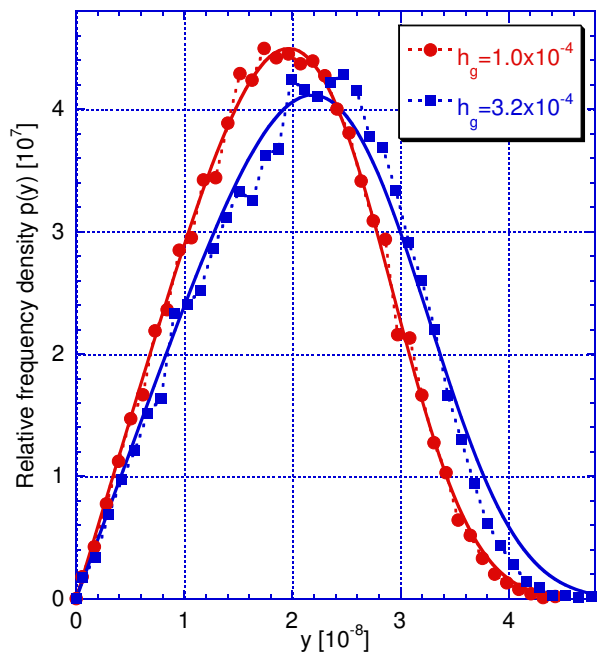,width=8.0cm}
  \caption{\textsl{Hitting distributions for VODE.}}
  \label{FigVODE}
\end{figure}
%
Possible internal parameters to be varied in solver packages more sophisticated than Euler's method are typically the
error tolerances \texttt{RTOL} (relative) and \texttt{ATOL} (absolute) and the global time step~$h_\text{g}$ (the time
interval after which the user requests solution output from the solver). However, naturally the user has no immediate
control over the size of the actual steps taken, which is determined algorithmically as a function of the error
tolerance parameters \texttt{RTOL} and \texttt{ATOL}, frequently by trial-and-error methods using heuristics, but only
rarely by an explicit formula. Nonetheless, as shown in Figure \ref{FigVODE}, distributions very similar to the ones
seen for Euler's algorithm (Figure \ref{FigEulerHitDistr}) can be generated. Experiments do not readily suggest a
simple relationship between the shape of the distribution (parameter~$a$ in equation \eqref{EqnFitDistr}) and any of
the parameters \texttt{ATOL}, \texttt{RTOL}, and~$h_\text{g}$. We suspect the lack of direct control over the time step
size to be the main reason for this behaviour. We found that in order to produce Figure \ref{FigVODE}, one has to use
$\texttt{RTOL}=0$, which we also attribute to the step adaptation.

For the solver RADAU5~\cite{HaiWan96} the results are qualitatively
similar, which supports the assertion that the observed phenomena
are not specific to a particular algorithm, but rather general effects.

\section{Conclusion}
We analysed the cumulative effect of rounding errors incurred by deterministic ODE solvers as the step size $h
\rightarrow 0$. We considered in particular the interesting case where the ordinary differential equation has a saddle
fixed point and showed that the numerical solution is inherently
random and also obtained its theoretical distribution in
terms of the time, step size and numerical precision. We showed that as the step size $h \rightarrow 0$, the numerical
solution exhibits three types of behaviour, depending on the time: initially it converges to the actual solution, it
then undergoes a transition stage, finally it diverges from the actual solution.

By performing multiple repetitions with different values of the time step size, we observed the random distributions
predicted theoretically. We demonstrated that during the transition period described above the numerical solution intersects
all the straight lines through the origin. The theoretical and numerical distributions for the points at
which these intersections occur showed very good agreement. Both the mean and the standard deviation of these
distributions were found to be of the form $a h^{\gamma}$, where $ \gamma \in (0, 1/2]$ is a constant determined by the
ODE system, and~$a$ was found explicitly in terms of the precision of the computer. We mainly focused on the explicit
Euler and RK4 methods, however, we also briefly considered the implicit solvers VODE and RADAU5 to demonstrate that the observed
effects are not specific to a particular numerical method.

\subsubsection*{Acknowledgments}

This work has been partially funded by the EPSRC (grant number GR/R85662/01) under the title ``Mathematical and
Numerical Analysis of Coagulation-Diffusion Processes in Chemical Engineering''. The authors thank James R. Norris and
Markus Kraft for suggesting the project and the collaboration, and helpful discussions.



\bibliographystyle{plainnat}

\vspace{2ex}

\noindent Sebastian Mosbach \\
Department of Chemical Engineering \\
University of Cambridge \\
Pembroke Street \\
Cambridge \\
CB2 3RA \\
UK \\
E-mail: sm453@cam.ac.uk

\noindent Amanda Turner \\
Statistical Laboratory \\ 
Centre for Mathematical Sciences \\
Wilberforce Road \\
Cambridge \\
CB3 0WB \\
UK \\
E-mail: A.G.Turner@statslab.cam.ac.uk 

\end{document}